\documentclass[11pt]{article}

\usepackage{latexsym}
\usepackage{amssymb}

\newtheorem{thm}{Theorem}[section]

\newtheorem{prop}{Proposition}[section]

\newtheorem{dfn}{Definition}[section]

\begin{document}
{

\begin{center}
\Large\bf
The two-dimensional moment problem in a strip.
\end{center}
\begin{center}
\bf Sergey M. Zagorodnyuk
\end{center}

\section{Introduction.}
In this paper we consider the following problem:
to find a non-negative Borel measure $\mu$ in a strip
$$ \Pi = \Pi(R) = \{ (x_1,x_2)\in \mathbb{R}^2:\ |x_2| \leq R \},\qquad R>0, $$
such that
\begin{equation}
\label{f1_1}
\int_\Pi x_1^m x_2^n d\mu = s_{m,n},\qquad m,n\in \mathbb{Z}_+,
\end{equation}
where $\{ s_{m,n} \}_{m,n\in \mathbb{Z}_+}$ is a prescribed sequence of complex numbers.
This problem is said to be {\bf the two-dimensional moment problem in a strip}.

The two-dimensional moment problem and the complex moment problem have
an extensive literature, see books~\cite{cit_10000_ST},
\cite{cit_450_Akh}, surveys~\cite{cit_455_F},\cite{cit_460_B} and~\cite{cit_470_S}.
However, to the best of our knowledge,
the two-dimensional moment problem in a strip was not solved.

\noindent
Firstly, we obtain a solvability criterion for the two-dimensional moment problem in a strip.
We  describe canonical solutions of this moment problem (see the definition
below). Secondly, we parameterize all solutions of the moment problem.
In a consequence, we derive conditions of the solvability and describe all solutions of
the complex moment problem with the support in a strip.
We shall use an abstract operator approach~\cite{cit_1500_Z} and results of Godi\v{c}, Lucenko and
Shtraus~\cite{cit_2000_GL},\cite[Theorem 1]{cit_3000_GP},\cite{cit_4000_S}.

{\bf Notations. } As usual, we denote by $\mathbb{R},\mathbb{C},\mathbb{N},\mathbb{Z},\mathbb{Z}_+$
the sets of real numbers, complex numbers, positive integers, integers and non-negative integers,
respectively.
For a subset $S$ of the complex plane we denote by $\mathfrak{B}(S)$ the set of all Borel subsets of $S$.
Everywhere in this paper, all Hilbert spaces are assumed to be separable. By
$(\cdot,\cdot)_H$ and $\| \cdot \|_H$ we denote the scalar product and the norm in a Hilbert space $H$,
respectively. The indices may be ommited in obvious cases.
For a set $M$ in $H$, by $\overline{M}$ we mean the closure of $M$ in the norm $\| \cdot \|_H$. For
$\{ x_{m,n} \}_{m,n\in \mathbb{Z}_+}$, $x_{m,n}\in H$, we write
$\mathop{\rm Lin}\nolimits \{ x_{m,n} \}_{m,n\in \mathbb{Z}_+}$ for the span of vectors
$\{ x_{m,n} \}_{m,n\in \mathbb{Z}_+}$
and $\mathop{\rm span}\nolimits \{ x_{m,n} \}_{m,n\in \mathbb{Z}_+} =
\overline{ \mathop{\rm Lin}\nolimits \{ x_{m,n} \}_{m,n\in \mathbb{Z}_+} }$.
The identity operator in $H$ is denoted by $E_H$. For an arbitrary linear operator $A$ in $H$,
the operators $A^*$,$\overline{A}$,$A^{-1}$ mean its adjoint operator, its closure and its inverse
(if they exist). By $D(A)$ and $R(A)$ we mean the domain and the range of the operator $A$.
By $\sigma(A)$, $\rho(A)$ we denote the spectrum of $A$ and the resolvent set of $A$, respectively.
We denote by $R_z (A)$ the resolvent function of $A$, $z\in \rho(A)$;
$\Delta_A(z) = (A-zE_H)D(A)$, $z\in \mathbb{C}$.
The norm of a bounded operator $A$ is denoted by $\| A \|$.
By $P^H_{H_1} = P_{H_1}$ we mean the operator of orthogonal projection in $H$ on a subspace
$H_1$ in $H$. By $\mathbf{B}(H)$ we denote the set of all bounded operators in $H$.

\section{Solvability of the two-dimensional moment problem in a strip.}
Let the two-dimensional moment problem in a strip~(\ref{f1_1}) be given.
Suppose that the moment problem has a solution $\mu$. Choose an arbitrary
polynomial $p(x_1,x_2)$ of the following form:
\begin{equation}
\label{f1_2}
\sum_{m=0}^\infty \sum_{n=0}^\infty \alpha_{m,n} x_1^m x_2^n,\qquad \alpha_{m,n}\in \mathbb{C},
\end{equation}
where all but finite number of coefficients $\alpha_{m,n}$ are zeros.
We may write
$$ 0 \leq \int_\Pi |p(x_1,x_2)|^2 d\mu =
\int_\Pi  \sum_{m,n=0}^\infty  \alpha_{m,n} x_1^m x_2^{n}
\overline{
\sum_{k,l=0}^\infty \alpha_{k,l} x_1^k x_2^{l}
} d\mu $$
$$ = \sum_{m,n,k,l} \alpha_{m,n}\overline{\alpha_{k,l}} \int_\Pi x_1^{m+k} x_2^{n+l} d\mu =
\sum_{m,n,k,l} \alpha_{m,n}\overline{\alpha_{k,l}} s_{m+k,n+l}. $$
Thus, for arbitrary complex numbers $\alpha_{m,n}$ (where all but finite numbers are zeros) we have
\begin{equation}
\label{f2_1}
\sum_{m,n,k,l=0}^\infty \alpha_{m,n}\overline{\alpha_{k,l}} s_{m+k,n+l} \geq 0.
\end{equation}
Since
$$ \int_\Pi |x_2 p(x_1,x_2)|^2 d\mu \leq R^2 \int_\Pi |p(x_1,x_2)|^2 d\mu, $$
in a similar manner we get
\begin{equation}
\label{f2_2}
\sum_{m,n,k,l=0}^\infty \alpha_{m,n}\overline{\alpha_{k,l}} (R^2 s_{m+k,n+l} - s_{m+k,n+l+2}) \geq 0,
\end{equation}
for arbitrary complex numbers $\alpha_{m,n}$ (where all but finite numbers are zeros).

On the other hand, suppose that the moment problem~(\ref{f1_1}) is given and conditions~(\ref{f2_1})
and~(\ref{f2_2}) hold. Let us show that the moment problem has a solution.
Set
\begin{equation}
\label{f2_2_1}
K((m,n),(k,l)) = s_{m+k,n+l},\qquad m,n,k,l\in \mathbb{Z}_+.
\end{equation}
Then relations~(\ref{f2_1}) may be written as
\begin{equation}
\label{f2_2_2}
\sum_{m,n,k,l=0}^\infty \alpha_{m,n}\overline{\alpha_{k,l}} K((m,n),(k,l)) \geq 0,
\end{equation}
for arbitrary complex numbers $\alpha_{m,n}$ (where all but finite numbers are zeros).
In this case $K$ is said  to be a {\it positive definite kernel on $\mathbb{Z}_+\times \mathbb{Z}_+$}.

We shall use the following important fact (e.g.~\cite[pp.361-363]{cit_6000_AG}).
\begin{thm}
\label{t2_1}
Let $K$ be a positive definite kernel on $\mathbb{Z}_+\times \mathbb{Z}_+$.
Then there exist a separable Hilbert space $H$ with a scalar product $(\cdot,\cdot)$ and
a sequence $\{ x_{m,n} \}_{m,n\in \mathbb{Z}_+}$ in $H$, such that
\begin{equation}
\label{f2_4}
K((m,n),(k,l)) = (x_{m,n},x_{k,l}),\qquad m,n,k,l\in \mathbb{Z}_+,
\end{equation}
and $\mathop{\rm span}\nolimits\{ x_{m,n} \}_{m,n\in \mathbb{Z}_+} = H$.
\end{thm}
{\bf Proof. }
Choose an arbitrary infinite-dimensional linear vector space $V$ (for instance, we may choose the space of all complex
sequences $(u_n)_{n\in \mathbb{N}}$, $u_n\in \mathbb{C}$).
Let $X = \{ x_{m,n} \}_{m,n\in \mathbb{Z}_+}$ be an arbitrary infinite sequence of linear independent elements
in $V$ which is indexed by elements  of $\mathbb{Z}_+\times \mathbb{Z}_+$.
Set $L_X = \mathop{\rm Lin}\nolimits\{ x_{m,n} \}_{m,n\in \mathbb{Z}_+}$. Introduce the following functional:
\begin{equation}
\label{f2_5}
[x,y] = \sum_{m,n,k,l=0}^\infty K((m,n),(k,l)) a_{m,n}\overline{b_{k,l}},
\end{equation}
for $x,y\in L_X$,
$$ x=\sum_{m,n=0}^\infty a_{m,n} x_{m,n},\quad y=\sum_{k,l=0}^\infty b_{k,l} x_{k,l},\quad
a_{m,n},b_{k,l}\in \mathbb{C}. $$
Here all but finite number of indices $a_{m,n},b_{k,l}$ are zeros.

\noindent
The set $L_X$ with $[\cdot,\cdot]$ will be a pre-Hilbert space. Factorizing and making the completion
we obtain the  desired space $H$ (\cite[p. 10-11]{cit_7000_B}).
$\Box$

By Theorem~\ref{t2_1} we obtain a Hilbert space $H$ and a sequence
$\{ x_{m,n} \}_{m,n\in \mathbb{Z}_+}$, $x_{m,n}\in H$, such that
\begin{equation}
\label{f2_6}
(x_{m,n}, x_{k,l})_H = K((m,n),(k,l)),\qquad m,n,k,l\in \mathbb{Z}_+.
\end{equation}
Set $L = \mathop{\rm Lin}\nolimits\{ x_{m,n} \}_{m,n\in \mathbb{Z}_+}$.
Introduce the following operators
\begin{equation}
\label{f2_7}
A_0 x = \sum_{m,n\in \mathbb{Z}_+} \alpha_{m,n} x_{m+1,n},
\end{equation}
\begin{equation}
\label{f2_8}
B_0 x = \sum_{m,n\in \mathbb{Z}_+} \alpha_{m,n} x_{m,n+1},
\end{equation}
where
\begin{equation}
\label{f2_9}
x = \sum_{m,n\in \mathbb{Z}_+} \alpha_{m,n} x_{m,n} \in L.
\end{equation}
Let us check that these definitions are correct.
Indeed, suppose that the element $x$ in~(\ref{f2_9}) has another representation:
\begin{equation}
\label{f2_10}
x = \sum_{k,l\in \mathbb{Z}_+} \beta_{k,l} x_{k,l}.
\end{equation}
We may write
$$ \left( \sum_{m,n\in \mathbb{Z}_+} \alpha_{m,n} x_{m+1,n}, x_{a,b} \right) =
\sum_{m,n\in \mathbb{Z}_+} \alpha_{m,n} K((m+1,n),(a,b)) $$
$$= \sum_{m,n\in \mathbb{Z}_+} \alpha_{m,n} s_{m+1+a,n+b} = \sum_{m,n\in \mathbb{Z}_+}
\alpha_{m,n} K((m,n),(a+1,b)) $$
$$ = \left(\sum_{m,n\in \mathbb{Z}_+} \alpha_{m,n} x_{m,n}, x_{a+1,b} \right) = (x,x_{a+1,b}), $$
for arbitrary $a,b\in \mathbb{Z}_+$.
In the same manner we get
$$ \left(\sum_{k,l\in \mathbb{Z}_+} \beta_{k,l} x_{k+1,l}, x_{a,b} \right) = (x,x_{a+1,b}). $$
Since $\mathop{\rm span}\nolimits\{ x_{a,b} \}_{a,b\in \mathbb{Z}_+} = H$, we get
$$ \sum_{m,n\in \mathbb{Z}_+} \alpha_{m,n} x_{m+1,n} =
\sum_{k,l\in \mathbb{Z}_+} \beta_{k,l} x_{k+1,l}. $$
Therefore the operator $A_0$ is defined correctly.
For $B_0$ considerations are similar.  It is not hard to see that operators $A_0$ and $B_0$ are
symmetric. Moreover, condition~(\ref{f2_2}) implies that the operator $B_0$ is bounded.
Set
\begin{equation}
\label{f2_10_1}
A = \overline{A_0},\quad B = \overline{B_0}.
\end{equation}
Observe that $B_0$ is a bounded self-adjoint operator in $H$. Since $A_0$ and $B_0$ commute, we easily get
\begin{equation}
\label{f2_11}
AB x = BA x,\qquad x\in D(A).
\end{equation}
We shall also need the following operator:
\begin{equation}
\label{f2_12}
J_0 x = \sum_{m,n\in \mathbb{Z}_+} \overline{\alpha_{m,n}} x_{m,n},
\end{equation}
where
\begin{equation}
\label{f2_13}
x = \sum_{m,n\in \mathbb{Z}_+} \alpha_{m,n} x_{m,n} \in L.
\end{equation}
Let us check that this definition is correct. Consider another representation for $x$ as in~(\ref{f2_10}).
Then
$$ \left\| \sum_{m,n\in \mathbb{Z}_+} (\overline{\alpha_{m,n}} - \overline{\beta_{m,n}}) x_{m,n} \right\|^2 $$
$$= \left( \sum_{m,n\in \mathbb{Z}_+} \overline{ (\alpha_{m,n}-\beta_{m,n}) } x_{m,n},
\sum_{k,l\in \mathbb{Z}_+} \overline{ (\alpha_{k,l}-\beta_{k,l}) } x_{k,l} \right) $$
$$ = \sum_{m,n,k,l\in \mathbb{Z}_+} \overline{(\alpha_{m,n}-\beta_{m,n})} (\alpha_{k,l}-\beta_{k,l})
K((m,n),(k,l)) $$
$$ = \sum_{m,n,k,l\in \mathbb{Z}_+} \overline{(\alpha_{m,n}-\beta_{m,n})} (\alpha_{k,l}-\beta_{k,l})
K((k,l),(m,n)) $$
$$= \left( \sum_{m,n\in \mathbb{Z}_+} (\alpha_{k,l}-\beta_{k,l}) x_{k,l},
\sum_{m,n\in \mathbb{Z}_+} (\alpha_{m,n}-\beta_{m,n}) x_{m,n} \right)  = 0. $$
Therefore the definition of $J_0$ is correct.
For an arbitrary $y = \sum_{a,b\in \mathbb{Z}_+} \gamma_{a,b} x_{a,b} \in L$ we may write
$$ (J_0 x,J_0 y) = \sum_{m,n,a,b} \overline{\alpha_{m,n}}\gamma_{a,b} (x_{m,n},x_{a,b}) =
\sum_{m,n,a,b} \overline{\alpha_{m,n}}\gamma_{a,b} K((m,n),(a,b)) $$
$$ = \sum_{m,n,a,b} \overline{\alpha_{m,n}} \gamma_{a,b} K((a,b),(m,n)) =
 \sum_{m,n,a,b} \overline{\alpha_{m,n}}\gamma_{a,b} (x_{a,b},x_{m,n}) =  (y,x). $$
In particular, this implies that $J_0$ is bounded. By continuity we extend $J_0$ to a bounded antilinear
operator $J$ such that
$$ (Jx,Jy) = (y,x),\qquad x,y\in H. $$
Moreover, we get $J^2 = E_H$. Consequently, $J$ is a conjugation in $H$ (\cite{cit_8000_S}).

\noindent
Notice that $J_0$ commutes with $A_0$ and $B_0$. It is easy to check that
\begin{equation}
\label{f2_14}
AJ x = JA x,\qquad x\in D(A),
\end{equation}
and by continuity we have
\begin{equation}
\label{f2_15}
BJ x = JB x,\qquad x\in H.
\end{equation}
Consider the Cayley transformations of the operators A,B:
\begin{equation}
\label{f2_16}
V_A := (A+iE_H)(A-iE_H)^{-1} = E + 2i(A-iE_H)^{-1},
\end{equation}
\begin{equation}
\label{f2_16_1}
U_B := (B+iE_H)(B-iE_H)^{-1} = E + 2i(B-iE_H)^{-1}.
\end{equation}
Set
\begin{equation}
\label{f2_17}
H_1 := \Delta_A(i),\ H_2 := H\ominus H_1,\ H_3:= \Delta_A(-i),\ H_4 := H\ominus H_3.
\end{equation}
\begin{prop}
\label{p2_1}
The operator $U_B$ reduces the subspaces $H_j$, $1\leq j\leq 4$:
\begin{equation}
\label{f2_18}
U_B H_j = H_j,\qquad 1\leq j\leq 4.
\end{equation}
Moreover, the following equality holds:
\begin{equation}
\label{f2_19}
U_B V_A x = V_A U_B x,\qquad x\in H_1.
\end{equation}
\end{prop}
{\bf Proof. } Choose an arbitrary $x\in \Delta_A(z)$, $x=(A-zE_H)f_A$, $f_A\in D(A)$,
$z\in \mathbb{C}\backslash \mathbb{R}$.
By~(\ref{f2_11}) we get
$$ Bx = BAf_A - zBf_A = ABf_A - zBf_A = (A-zE_H)Bf_A\in \Delta_A(z). $$
In particular, we have
\begin{equation}
\label{f2_19_1}
BH_1\subseteq H_1,\quad BH_3\subseteq H_3.
\end{equation}
Since $B$ is bounded and self-adjoint, we get
\begin{equation}
\label{f2_19_2}
BH_2\subseteq H_2,\quad BH_4\subseteq H_4.
\end{equation}
This implies equalities
\begin{equation}
\label{f2_19_3}
(B-iE_H) H_j = H_j,\quad 1\leq j\leq 4;
\end{equation}
\begin{equation}
\label{f2_19_4}
(B-iE_H)^{-1} H_j = H_j,\quad 1\leq j\leq 4.
\end{equation}
Then $U_B H_j\subseteq H_j$, $1\leq j\leq 4$, and relation~(\ref{f2_18}) follows.

\noindent
Since
$$ (A-iE_H) Bx = B(A-iE_H)x,\quad x\in D(A),$$
we get
\begin{equation}
\label{f2_19_5}
BV_A y = V_A B y,\quad y\in H_1.
\end{equation}
Then
$$ (B-iE_H) V_A y = V_A (B-iE_H) y,\quad y\in H_1, $$
and, hence, we easily get relation~(\ref{f2_19}).
$\Box$

We shall construct a unitary operator $U$ in $H$, $U\supset V_A$, which commutes with $U_B$.
Choose an arbitrary $x\in H$, $x= x_{H_1} + x_{H_2}$, $x_{H_1}\in H_1$, $x_{H_2}\in H_2$.
If $U$ is an operator of the required type then
(we use Proposition~\ref{p2_1}):
$$ U_B U x = U_B V_A x_{H_1} + U_B U x_{H_2} = V_A U_B x_{H_1} + U_B U x_{H_2}, $$
$$ U U_B x = U U_B x_{H_1} + U U_B x_{H_2} = V_A U_B x_{H_1} + U U_B x_{H_2}. $$
Thus, we have to find an isometric operator $U_{2,4}$ which maps $H_2$ onto $H_4$, and
commutes with $U_B$:
\begin{equation}
\label{f2_21}
U_B U_{2,4} x = U_{2,4} U_B x,\qquad x\in H_2.
\end{equation}
Moreover, all operators $U$ of the required type have the following form:
\begin{equation}
\label{f2_22}
U = V_A \oplus U_{2,4},
\end{equation}
where $U_{2,4}$ is an isometric operator which maps $H_2$ onto $H_4$, and
commutes with $U_B$.

\noindent
Denote the operator $U_B$ restricted to $H_i$ by $U_{B;H_i}$, $1\leq i\leq 4$.
Notice that
\begin{equation}
\label{f2_23}
A^* J x= JA^* x,\qquad x\in D(A^*).
\end{equation}
Indeed, for arbitrary $f_A\in D(A)$ and $g_{A^*}\in D(A^*)$ we may write
$$ \overline{ (Af_A,Jg_{A^*})  } = (JAf_A, g_{A^*}) = (AJf_A, g_{A^*}) = (Jf_A, A^*g_{A^*}) $$
$$ = \overline{ (f_A,JA^*g_{A^*})  }, $$
and~(\ref{f2_23}) follows.

\noindent
Choose an arbitrary $x\in H_2$. We have
$$ A^* x = -i x, $$
and therefore
$$ A^* Jx = JA^* x = ix. $$
Thus, we have
$$ JH_2 \subseteq H_4. $$
In a similar manner we get
$$ JH_4 \subseteq H_2, $$
and therefore
\begin{equation}
\label{f2_24}
JH_2 = H_4,\quad JH_4 = H_2.
\end{equation}
By the Godi\v{c}-Lucenko Theorem (\cite{cit_2000_GL},\cite[Theorem 1]{cit_3000_GP}) we have a
representation:
\begin{equation}
\label{f2_25}
U_{B;H_2} = KL,
\end{equation}
where $K$ and $L$ are some conjugations in $H_2$.
We set
\begin{equation}
\label{f2_26}
U_{2,4} := JK.
\end{equation}
From~(\ref{f2_24}) it follows that $U_{2,4}$ maps isometrically $H_2$ onto $H_4$.
Observe that
\begin{equation}
\label{f2_27}
U_{2,4}^{-1} := KJ.
\end{equation}
Notice that
\begin{equation}
\label{f2_27_1}
JU_BJ = U_B^{-1}.
\end{equation}
Indeed, by virtue of~\cite[Proposition 2.10]{cit_10000_Z} we can write
$$ JU_BJ = E - 2iJ(B-iE_H)^{-1}J = E - 2i(J(B-iE_H)J)^{-1} $$
$$ = E- 2i(B+iE_H)^{-1} = E - 2iR_B(-i) = U_B^* = U_B^{-1}. $$
By~(\ref{f2_27_1}) we get
$$ U_{2,4} U_{B;H_2} U_{2,4}^{-1} x = JK KL KJ x = J LK J x = J U_{B;H_2}^{-1} J x $$
$$ = J U_B^{-1}J x = U_B x = U_{B;H_4} x,\qquad x\in H_4. $$
Therefore relation~(\ref{f2_21}) is true.
We define an operator $U$ by~(\ref{f2_22}) and set
\begin{equation}
\label{f2_28}
A_U := i(U+E_H)(U-E_H)^{-1} = iE_H + 2i(U-E_H)^{-1}.
\end{equation}
The inverse Cayley transformation $A_U$ is correctly defined since $1$ is not in the point spectrum of $U$.
Indeed, $V_A$ is the Cayley transformation of a symmetric operator while eigen subspaces $H_2$ and
$H_4$ have the zero intersection.
Let
\begin{equation}
\label{f2_29}
A_U = \int_\mathbb{R} x_1 dE(x_1),\quad B = \int_{ [-R,R] } x_2 dF(x_2),
\end{equation}
where $E$ and $F$ are the spectral measures of $A_U$ and $B$, respectively. These measures are
defined on $\mathfrak{B}(\mathbb{R})$ (\cite{cit_9000_BS}).
Since $U$ and $U_B$ commute, we get that $E$ and $F$ commute, as well.
By the induction argument we get
$$ x_{m,n} = A^m x_{0,n},\qquad m,n\in \mathbb{Z}_+, $$
and
$$ x_{0,n} = B^n x_{0,0},\qquad n\in \mathbb{Z_+}. $$
Therefore we obtain
\begin{equation}
\label{f2_30}
x_{m,n} = A^m B^n x_{0,0},\qquad m,n\in \mathbb{Z}_+.
\end{equation}
We may write
$$ x_{m,n} = \int_\mathbb{R} x_1^m dE(x_1) \int_{ [-R,R] } x_2^n dF(x_2) x_{0,0} =
\int_\Pi x_1^m x_2^n  d(E\times F)(x_1,x_2) x_{0,0}, $$
where $E\times F$ is the product spectral measure on $\mathfrak{B}(\Pi)$.
Then
\begin{equation}
\label{f2_31}
s_{m,n} = (x_{m,n},x_{0,0})_H = \int_\Pi x_1^m x_2^n d((E\times F) x_{0,0}, x_{0,0})_H,\quad
m,n\in \mathbb{Z}_+.
\end{equation}
The measure $\mu := ((E\times F) x_{0,0}, x_{0,0})_H$ is a non-negative Borel measure on $\Pi$ and
relation~(\ref{f2_31}) shows that $\mu$ is a solution of the moment problem~(\ref{f1_1}).

\begin{thm}
\label{t2_2}
Let the moment problem~(\ref{f1_1}) be given.
This problem has a solution if an only if conditions~(\ref{f2_1}),(\ref{f2_2}) hold
for arbitrary complex numbers $\alpha_{m,n}$ such that all but finite numbers are zeros.
\end{thm}

\section{A parameterization of all solutions of the two-dimensional moment problem in a strip.}
Let the moment problem~(\ref{f1_1}) be given. Define a Hilbert space $H$ and operators
$A,B,J$ as in the previous Section.
Let $\widetilde A\supseteq A$ be a self-adjoint extension of $A$ in a Hilbert space
$\widetilde H\supseteq H$ and $E_{\widetilde A}$ be the spectral measure of $\widetilde A$.
Recall that the function
\begin{equation}
\label{f3_1}
\mathbf{R}_z(A) := P^{\widetilde H}_H R_z(\widetilde A),\qquad z\in \mathbb{C}\backslash \mathbb{R},
\end{equation}
is said to be a {\it generalized resolvent} of $A$. The function
\begin{equation}
\label{f3_2}
\mathbf{E}_A (\delta) := P^{\widetilde H}_H E_{\widetilde A} (\delta),\qquad \delta\in \mathfrak{B}(\mathbb{R}),
\end{equation}
is said to be a {\it spectral measure}  of $A$.
There exists a one-to-one correspondence between generalized resolvents and spectral measures according
to the following relation~\cite{cit_6000_AG}:
\begin{equation}
\label{f3_3}
(\mathbf{R}_z(A) x,y)_H = \int_{\mathbb{R}} \frac{1}{t-z} d(\mathbf{E}_A x,y)_H,\qquad x,y\in H.
\end{equation}

\begin{thm}
\label{t3_1}
Let the moment problem~(\ref{f1_1}) be given and conditions~(\ref{f2_1}),(\ref{f2_2}) hold.
Consider a Hilbert space $H$ and a sequence
$\{ x_{m,n} \}_{m,n\in \mathbb{Z}_+}$, $x_{m,n}\in H$, such that relation~(\ref{f2_6})
holds where $K$ is defined by~(\ref{f2_2}).
Consider operators $A_0$,$B_0$,$A$,$B$ defined by~(\ref{f2_7}),(\ref{f2_8}) and~(\ref{f2_10_1}).
Let $\mu$ be an arbitrary solution of the moment problem. Then  it has the following form:
\begin{equation}
\label{f3_4}
\mu (\delta)= ((\mathbf{E}\times F)(\delta) x_{0,0}, x_{0,0})_H,\qquad \delta\in \mathfrak{B}(\Pi),
\end{equation}
where $F$ is the spectral measure of $B$, $\mathbf{E}$ is a spectral measure of $A$ which commutes with
$F$. By $((\mathbf{E}\times F)(\delta) x_{0,0}, x_{0,0})_H$ we mean the non-negative Borel measure on
$\mathbb{R}$ which is obtained by the Lebesgue continuation procedure from the following
non-negative measure on rectangules
\begin{equation}
\label{f3_5}
((\mathbf{E}\times F)(I_{x_1}\times I_{x_2}) x_{0,0}, x_{0,0})_H :=
( \mathbf{E}(I_{x_1}) F(I_{x_2}) x_{0,0}, x_{0,0})_H,
\end{equation}
where $I_{x_1}\subset \mathbb{R}$, $I_{x_2}\subseteq [-R,R]$ are arbitrary intervals.

\noindent
On the other hand, for an arbitrary spectral measure $\mathbf{E}$ of $A$ which commutes with the
spectral measure $F$ of $B$, by relation~(\ref{f3_4}) it corresponds a solution of the moment
problem~(\ref{f1_1}).

\noindent
Moreover, the correspondence between the spectral measures of $A$ which commute with the spectral meeasure of
$B$ and solutions of the moment problem is bijective.
\end{thm}
{\bf Remark. } It is straightforward to check that the measure in~(\ref{f3_5}) is non-negative and additive.
Moreover, the standard arguments~\cite[Chapter 5, Theorem 2, p. 254-255]{cit_9500_KF} imply
that the measure in~(\ref{f3_5}) is $\sigma$-additive.
Consequently, it has the (unique) Lebesgue continuation to a (finite) non-negative Borel measure
on $\Pi$.

{\bf Proof. }
Consider a Hilbert space $H$ and operators $A$,$B$ as in the statement of the Theorem.
Let $F$ be the spectral measure of $B$. Let $\mu$ be an arbitrary solution of the moment problem~(\ref{f1_1}).
Consider the space $L^2_\mu$ of complex functions on $\Pi$ which are square integrable with respect to the
measure $\mu$. The scalar product and the norm are given by
$$ (f,g)_\mu =
\int_\Pi f(x_1,x_2) \overline{ g(x_1,x_2) } d\mu,\quad
\|f\|_\mu = \left( (f,f)_\mu \right)^{ \frac{1}{2} },\quad f,g\in L^2_\mu. $$
Consider the following operators:
\begin{equation}
\label{f3_6}
A_\mu f(x_1,x_2) = x_1 f(x_1,x_2),\qquad D(A_\mu) = \{ f\in L^2_\mu:\ x_1 f(x_1,x_2)\in L^2_\mu \},
\end{equation}
\begin{equation}
\label{f3_7}
B_\mu f(x_1,x_2) = x_2 f(x_1,x_2),\qquad D(B_\mu) = L^2_\mu.
\end{equation}
The operator $A_\mu$ is self-adjoint and the operator $B_\mu$ is self-adjoint and bounded.
These operators commute and therefore the spectral measure $
E_\mu$ of $A_\mu$ and the spectral measure $F_\mu$ of $B_\mu$
commute, as well.

\noindent
Let $p(x_1,x_2)$ be a  polynomial of the form~(\ref{f1_1}) and
$q(x_1,x_2)$ be a polynomial of the form~(\ref{f1_1}) with
$\beta_{m,n}\in \mathbb{C}$ instead of $\alpha_{m,n}$.
Then
$$ (p,q)_\mu = \sum_{m,n,k,l\in \mathbb{Z}_+} \alpha_{m,n}\overline{ \beta_{k,l} }
\int_\Pi x_1^{m+k} x_2^{n+l} d\mu $$
$$ = \sum_{m,n,k,l\in \mathbb{Z}_+} \alpha_{m,n}\overline{ \beta_{k,l} } s_{m+k,n+l}, $$
On the other hand, we may write
$$ \left(
\sum_{m,n\in \mathbb{Z}_+} \alpha_{m,n} x_{m,n}, \sum_{k,l\in \mathbb{Z}_+} \beta_{k,l} x_{k,l} \right)_H =
\sum_{m,n,k,l\in \mathbb{Z}_+} \alpha_{m,n}\overline{ \beta_{k,l} }
(x_{m,n},x_{k,l})_H $$
$$ = \sum_{m,n,k,l\in \mathbb{Z}_+} \alpha_{m,n}\overline{ \beta_{k,l} } K((m,n),(k,l))
= \sum_{m,n,k,l\in \mathbb{Z}_+} \alpha_{m,n}\overline{ \beta_{k,l} } s_{m+k,n+l}. $$
Therefore
\begin{equation}
\label{f3_8}
(p,q)_\mu = \left(
\sum_{m,n\in \mathbb{Z}_+} \alpha_{m,n} x_{m,n}, \sum_{k,l\in \mathbb{Z}_+} \beta_{k,l} x_{k,l} \right)_H.
\end{equation}
Consider the following operator:
\begin{equation}
\label{f3_9}
V[p] = \sum_{m,n\in \mathbb{Z}_+} \alpha_{m,n} x_{m,n},\quad p=\sum_{m,n\in \mathbb{Z}_+}
\alpha_{m,n} x_1^m x_2^n.
\end{equation}
Here by $[p]$ we mean the class of equivalence in $L^2_\mu$ defined by $p$. If two different polynomials
$p$ and $q$ belong to the same class of equivalence then by~(\ref{f3_8}) we get
$$ 0 = \| p-q \|_\mu^2 = (p-q,p-q)_\mu = \left( \sum_{m,n\in \mathbb{Z}_+} (\alpha_{m,n}-\beta_{m,n}) x_{m,n},
\sum_{k,l\in \mathbb{Z}_+} (\alpha_{k,l}-\beta_{k,l}) x_{k,l} \right)_H $$
$$ = \left\| \sum_{m,n\in \mathbb{Z}_+} \alpha_{m,n} x_{m,n} -
\sum_{m,n\in \mathbb{Z}_+} \beta_{m,n} x_{m,n} \right\|_H^2. $$
Thus, the definition of $V$ is correct.
The operator $V$ maps the set of all polynomials $P^2_{0,\mu}$
in $L^2_\mu$ on $L$. By continuity we extend $V$ to an isometric transformation from
the closure of polynomials $P^2_\mu = \overline{P^2_{0,\mu}}$ onto $H$.

\noindent
Set $H_0 := L^2_\mu \ominus P^2_\mu$. Introduce the following operator:
\begin{equation}
\label{f3_10}
U := V \oplus E_{H_0},
\end{equation}
which maps isometrically $L^2_\mu$ onto $\widetilde H := H\oplus H_0$.
Set
\begin{equation}
\label{f3_11}
\widetilde A := UA_\mu U^{-1},\quad \widetilde B := UB_\mu U^{-1}.
\end{equation}
Notice that
$$ \widetilde A x_{m,n} = UA_\mu U^{-1} x_{m,n} = UA_\mu x_1^m x_2^n = Ux_1^{m+1} x_2^{n} = x_{m+1,n}, $$
$$ \widetilde B x_{m,n} = UB_\mu U^{-1} x_{m,n} = UB_\mu x_1^m x_2^n = Ux_1^{m} x_2^{n+1} = x_{m,n+1}. $$
Therefore $\widetilde A\supseteq A$ and $\widetilde B\supseteq B$.
Let
\begin{equation}
\label{f3_12}
\widetilde A = \int_\mathbb{R} x_1 d\widetilde E(x_1),\quad
\widetilde B = \int_{ [-R,R] } x_2
d \widetilde F(x_2),
\end{equation}
where $\widetilde E$ and $\widetilde F$ are the spectral measures of $\widetilde A$ and
$\widetilde B$, respectively.
Repeating arguments after relation~(\ref{f2_29}) we obtain that
\begin{equation}
\label{f3_13}
x_{m,n} = \widetilde A^m \widetilde B^n x_{0,0},\qquad m,n\in \mathbb{Z}_+,
\end{equation}
\begin{equation}
\label{f3_14}
s_{m,n} = \int_\Pi x_1^m x_2^n d((\widetilde E\times \widetilde F)(x_1,x_2) x_{0,0}, x_{0,0})_{\widetilde H},\quad
m,n\in \mathbb{Z}_+,
\end{equation}
where $(\widetilde E\times \widetilde F)$ is the product measure of $\widetilde E$ and $\widetilde F$.
Thus, the measure $\widetilde \mu := ((\widetilde E\times \widetilde F) x_{0,0}, x_{0,0})_{\widetilde H}$
is a solution of the moment problem.

\noindent
Let $I_{x_1}\subset \mathbb{R}$, $I_{x_2}\subseteq [-R,R]$ be arbitrary intervals.
Observe  that
$$ P^{\widetilde H}_H \widetilde E(I_{x_1}) \widetilde F(I_{x_2}) P^{\widetilde H}_H =
P^{\widetilde H}_H \widetilde E(I_{x_1}) P^{\widetilde H}_H \widetilde F(I_{x_2}) P^{\widetilde H}_H
= \mathbf{E}(I_{x_1}) F(I_{x_2}); $$
$$ P^{\widetilde H}_H \widetilde E(I_{x_1}) \widetilde F(I_{x_2}) P^{\widetilde H}_H =
P^{\widetilde H}_H \widetilde F(I_{x_2}) \widetilde E(I_{x_1}) P^{\widetilde H}_H =
P^{\widetilde H}_H \widetilde F(I_{x_2}) P^{\widetilde H}_H \widetilde E(I_{x_1}) P^{\widetilde H}_H $$
$$ = F(I_{x_2}) \mathbf{E}(I_{x_1}), $$
and therefore
\begin{equation}
\label{f3_14_1}
\mathbf{E}(I_{x_1}) F(I_{x_2}) = F(I_{x_2}) \mathbf{E}(I_{x_1}),
\end{equation}
where $\mathbf{E}$ is the corresponding spectral function of $A$ and $F$ is the spectral function of $B$.
Then
$$ \widetilde \mu (I_{x_1} \times I_{x_2}) = ((\widetilde E\times \widetilde F) (I_{x_1} \times I_{x_2})
x_{0,0}, x_{0,0})_{\widetilde H} $$
$$ = ( \widetilde E(I_{x_1}) \widetilde F(I_{x_2}) x_{0,0}, x_{0,0})_{\widetilde H} =
( P^{\widetilde H}_H \widetilde F(I_{x_2}) \widetilde E(I_{x_1}) x_{0,0}, x_{0,0})_{\widetilde H} $$
$$ = ( P^{\widetilde H}_H \widetilde F(I_{x_2}) P^{\widetilde H}_H
\widetilde E(I_{x_1}) x_{0,0}, x_{0,0})_{\widetilde H}
= ( F(I_{x_2}) \mathbf{E}(I_{x_1}) x_{0,0}, x_{0,0})_{H} $$
$$ = ( \mathbf{E}(I_{x_1}) F(I_{x_2}) x_{0,0}, x_{0,0})_{H}. $$
where $\mathbf{E}$ is the corresponding spectral function of $A$ and $F$ is the spectral function of $B$.
Thus, the measure $\widetilde \mu$ admits the representation~(\ref{f3_4}) since the Lebesgue continuation
is unique.

\noindent
Let us show that $\widetilde \mu = \mu$.
Consider the following transformation:
\begin{equation}
\label{f3_15}
S:\ (x_1,x_2) \in \Pi \mapsto \left( \mathop{\rm Arg }\nolimits \frac{x_1-i}{x_1+i}, x_2  \right) \in \Pi_0,
\end{equation}
where $\Pi_0 = [-\pi,\pi) \times [-R,R]$ and $\mathop{\rm Arg }\nolimits e^{iy} = y\in [-\pi,\pi)$.
By virtue of $V$ we define  the following measures:
\begin{equation}
\label{f3_16}
\mu_0 (VG) := \mu (G),\quad \widetilde\mu_0 (VG) := \widetilde\mu (G),\qquad G\in \mathfrak{B}(\Pi),
\end{equation}
It is not hard to see that $\mu_0$ and $\widetilde\mu_0$ are non-negative measures on
$\mathfrak{B}(\Pi_0)$.
Then
\begin{equation}
\label{f3_17}
\int_\Pi \left( \frac{x_1-i}{x_1+i} \right)^m x_2^n d\mu =
\int_{\Pi_0} e^{im\psi} x_2^n d\mu_0,
\end{equation}
\begin{equation}
\label{f3_18}
\int_\Pi \left( \frac{x_1-i}{x_1+i} \right)^m x_2^n d\widetilde\mu =
\int_{\Pi_0} e^{im\psi} x_2^n d\widetilde\mu_0,\qquad m\in \mathbb{Z},n\in \mathbb{Z}_+;
\end{equation}
and
$$  \int_\Pi \left( \frac{x_1-i}{x_1+i} \right)^m x_2^n d\widetilde\mu =
\int_\Pi \left( \frac{x_1-i}{x_1+i} \right)^m x_2^n
d((\widetilde E\times \widetilde F) x_{0,0}, x_{0,0})_{\widetilde H} $$
$$ = \left( \int_\Pi \left( \frac{x_1-i}{x_1+i} \right)^m x_2^n
d(\widetilde E\times \widetilde F) x_{0,0}, x_{0,0} \right)_{\widetilde H}  $$
$$ = \left( \int_\mathbb{R} \left( \frac{x_1-i}{x_1+i} \right)^m d\widetilde E
\int_{[-R,R]} x_2^n d\widetilde F x_{0,0}, x_{0,0} \right)_{\widetilde H}  $$
$$ = \left( \left( (\widetilde A - iE_{\widetilde H})(\widetilde A + iE_{\widetilde H})^{-1} \right)^m
\widetilde B^n  x_{0,0}, x_{0,0} \right)_{\widetilde H}  $$
$$ = \left( U^{-1}\left( (\widetilde A - iE_{\widetilde H})(\widetilde A + iE_{\widetilde H})^{-1} \right)^m
\widetilde B^n  U 1, U 1 \right)_\mu  $$
$$ = \left( \left( (A_\mu - iE_{L^2_\mu})(A_\mu + iE_{L^2_\mu})^{-1} \right)^m
B_\mu^n  1, 1 \right)_\mu  $$
\begin{equation}
\label{f3_19}
 = \int_\Pi \left( \frac{x_1-i}{x_1+i} \right)^m x_2^n d\mu,\qquad m\in \mathbb{Z},n\in \mathbb{Z}_+.
\end{equation}
By virtue of relations~(\ref{f3_17}),(\ref{f3_18}) and~(\ref{f3_19}) we get
\begin{equation}
\label{f3_20}
\int_{\Pi_0} e^{im\psi} x_2^{n} d\mu_0 =
\int_{\Pi_0} e^{im\psi} x_2^{n} d\widetilde\mu_0,\qquad m\in \mathbb{Z},n\in \mathbb{Z}_+.
\end{equation}
By the Weierstrass theorem we can approximate any continuous function by exponentials and therefore
\begin{equation}
\label{f3_21}
\int_{\Pi_0} f(\psi) x_2^{n} d\mu_0 =
\int_{\Pi_0} f(\psi) x_2^{n} d\widetilde\mu_0,\qquad n\in \mathbb{Z}_+,
\end{equation}
for arbitrary continuous functions on $\Pi_0$. In particular, we have
\begin{equation}
\label{f3_22}
\int_{\Pi_0} x_1^{m} x_2^{n} d\mu_0 =
\int_{\Pi_0} x_1^{m} x_2^{n} d\widetilde\mu_0,\qquad n,m\in \mathbb{Z}_+.
\end{equation}
However, the two-dimensional Hausdorff moment problem is determinate (\cite{cit_10000_ST}) and therefore we get
$\mu_0 = \widetilde\mu_0$ and $\mu=\widetilde\mu$.
Thus, we have proved that an arbitrary solution $\mu$ of the moment problem~(\ref{f1_1}) can be represented
in the form~(\ref{f3_4}).

Let us check the second assertion of the Theorem.
For an arbitrary spectral measure $\mathbf{E}$ of $A$ which commutes with the
spectral measure $F$ of $B$, by relation~(\ref{f3_4}) we define a non-negative Borel measure $\mu$
on $\Pi$. Let us show that the measure $\mu$ is a solution of the moment
problem~(\ref{f1_1}).

\noindent
Let $\widehat A$ be a self-adjoint extension of the operator $A$ in a Hilbert space
$\widehat H\supseteq H$, such that
$$ \mathbf{E} = P^{\widehat H}_H \widehat E, $$
where $\widehat E$ is the spectral measure of $\widehat A$.
By~(\ref{f2_30}) we get
$$ x_{m,n} = A^m B^n x_{0,0} = \widehat A^m B^n x_{0,0} = P^{\widehat H}_H \widehat A^m B^n x_{0,0} $$
$$ = P^{\widehat H}_H \lim_{a\to +\infty} \int_{[-a,a)} x_1^m d\widehat E(x_1)
B^n x_{0,0}
= \lim_{a\to +\infty} P^{\widehat H}_H \int_{[-a,a)} x_1^m d\widehat{E}(x_1)  $$
$$ * B^n x_{0,0}
= \lim_{a\to +\infty} \int_{[-a,a)} x_1^m d\mathbf{E} (x_1) B^n x_{0,0}, $$
\begin{equation}
\label{f3_23}
\qquad m,n\in \mathbb{Z}_+,
\end{equation}
where the integrals are understood as strong limits of the Stieltjes operator sums.
We choose arbitrary points
$$  -a = x_{1,0} < x_{1,1} < ... < x_{1,N}=a; $$
\begin{equation}
\label{f3_24}
\max_{1\leq i\leq N}|x_{1,i}-x_{1,i-1}| =: d,\quad  N\in \mathbb{N};
\end{equation}
$$  -R = x_{2,0} < x_{2,1} < ... < x_{2,M} = R; $$
\begin{equation}
\label{f3_25}
 \max_{1\leq j\leq M} |x_{2,j}-x_{2,j-1}| =: r;\quad M\in \mathbb{N}.
\end{equation}
Set
$I_{2,j} = [x_{2,j-1},x_{2,j})$, if $1\leq j< M$, and $I_{2,M} = [x_{2,M-1},x_{2,M}]$.
Then
$$ C_a := \int_{[-a,a)} x_1^m d\mathbf{E} \int_{[-R,R]} x_2^n dF =
\lim_{d\rightarrow 0} \sum_{i=1}^N x_{1,i-1}^m \mathbf{E}([x_{1,i-1},x_{1,i})) $$
$$ * \lim_{r\rightarrow 0} \sum_{j=1}^M x_{2,j-1}^n F(I_{2,j}), $$
where the integral sums converge in the strong operator topology. Then
$$ C_a = \lim_{d\rightarrow 0} \lim_{r\rightarrow 0} \sum_{i=1}^N x_{1,i-1}^m \mathbf{E}([x_{1,i-1},x_{1,i}))
\sum_{j=1}^M x_{2,j-1}^n F(I_{2,j}) $$
$$ = \lim_{d\rightarrow 0} \lim_{r\rightarrow 0}
\sum_{i=1}^N \sum_{j=1}^M
x_{1,i-1}^m x_{2,j-1}^n
\mathbf{E}([x_{1,i-1},x_{1,i})) F(I_{2,j}), $$
where the limits are understood in the strong operator topology. Then
$$ (C_a x_{0,0}, x_{0,0})_H =
\left( \lim_{d\rightarrow 0} \lim_{r\rightarrow 0}
\sum_{i=1}^N \sum_{j=1}^M
x_{1,i-1}^m x_{2,j-1}^n
\mathbf{E}([x_{1,i-1},x_{1,i})) F(I_{2,j}) x_{0,0}, x_{0,0} \right)_H $$
$$ = \lim_{d\rightarrow 0} \lim_{r\rightarrow 0}
\sum_{i=1}^N \sum_{j=1}^M
x_{1,i-1}^m x_{2,j-1}^n
\left( \mathbf{E}([x_{1,i-1},x_{1,i})) F(I_{2,j}) x_{0,0}, x_{0,0} \right)_H $$
$$ = \lim_{d\rightarrow 0} \lim_{r\rightarrow 0}
\sum_{i=1}^N \sum_{j=1}^M
x_{1,i-1}^m x_{2,j-1}^n
\left( (\mathbf{E}\times F) ( [x_{1,i-1},x_{1,i})\times I_{2,j} ) x_{0,0}, x_{0,0} \right)_H $$
$$ = \lim_{d\rightarrow 0} \lim_{r\rightarrow 0}
\sum_{i=1}^N \sum_{j=1}^M
x_{1,i-1}^m x_{2,j-1}^n
\left( \mu ( [x_{1,i-1},x_{1,i})\times I_{2,j}) x_{0,0}, x_{0,0} \right)_H. $$
Therefore
$$ (C_a x_{0,0}, x_{0,0})_H =
\lim_{d\rightarrow 0} \lim_{r\rightarrow 0}
\int_{[-a,a)\times[-R,R]} f_{d,r} (x_1,x_2) d\mu, $$
where $f_{d,r}$ is equal to $x_{1,i-1}^m x_{2,j-1}^n$ on the rectangular
$[x_{1,i-1},x_{1,i}) \times I_{2,j}$, $1\leq i\leq N$, $1\leq j\leq M$.

\noindent
If $r\rightarrow 0$, then the function
$f_{d,r}(x_1,x_2)$ converges pointwise to a function $f_d(x_1,x_2)$ which is equal to
$x_{1,i-1}^m x_2^n$ on the rectangular
$[x_{1,i-1},x_{1,i}) \times [-R,R]$, $1\leq i\leq N$. Moreover, the function $f_{d,r}(x_1,x_2)$
is uniformly bounded.
By the Lebesgue we obtain
$$ (C_a x_{0,0}, x_{0,0})_H =
\lim_{d\rightarrow 0}
\int_{[-a,a)\times[-R,R]} f_{d} (x_1,x_2) d\mu. $$
If $d\rightarrow 0$, then the function
$f_{d}$ converges pointwise to a function $x_1^m x_2^n$. Since
$|f_d|\leq a^m R^n$, by the Lebesgue theorem we get
\begin{equation}
\label{f3_26}
(C_a x_{0,0}, x_{0,0})_H =
\int_{[-a,a)\times[-R,R]} x_1^m x_2^n d\mu.
\end{equation}
By virtue of relations~(\ref{f3_23}) and~(\ref{f3_26}) we get
$$ s_{m,n} = (x_{m,n},x_{0,0})_H = \lim_{a\to +\infty} (C_a x_{0,0},x_{0,0})_H $$
\begin{equation}
\label{f3_27}
= \lim_{a\to+\infty} \int_{[-a,a)\times[-R,R]} x_1^m x_2^n d\mu =
\int_\Pi x_1^m x_2^n d\mu.
\end{equation}
Thus, the measure $\mu$ is a solution of the moment problem (\ref{f1_1}).

Let us prove the last assertion of the Theorem. Suppose to the contrary that two different
spectral measures $\mathbf{E}_1$ and $\mathbf{E}_1$ of $A$ commute with the spectral measure $F$ of
$B$ and produce by relation~(\ref{f3_4}) the same solution $\mu$ of the moment problem.
Choose an arbitrary $z\in \mathbb{C}\backslash \mathbb{R}$. Then
$$ \int_\Pi \frac{x_1^m}{x_1-z} x_2^n d\mu =
\int_\Pi \frac{x_1^m}{x_1-z} x_2^n ((\mathbf{E}_k\times F)(\delta) x_{0,0}, x_{0,0})_H $$
\begin{equation}
\label{f3_28}
= \lim_{a\to +\infty}
\int_{[-a,a)\times [-R,R]}
\frac{x_1^m}{x_1-z} x_2^n d((\mathbf{E}_k\times F)(\delta) x_{0,0}, x_{0,0})_H,\quad k=1,2.
\end{equation}
Consider arbitrary partitions of the type~(\ref{f3_24}),(\ref{f3_25}). Then
$$ D_a := \int_{[-a,a)\times [-R,R]}
\frac{x_1^m}{x_1-z} x_2^n d((\mathbf{E}_k\times F)(\delta) x_{0,0}, x_{0,0})_H $$
$$ = \lim_{d\to 0} \lim_{r\to 0}
\int_{[-a,a)\times [-R,R]} g_{z;d,r}(x,\varphi)
d((\mathbf{E}_k\times F)(\delta) x_{0,0}, x_{0,0})_H. $$
Here the function $g_{z;d,r}(x_1,x_2)$ is equal to
$\frac{ x_{1,i-1}^m }{ x_{1,i-1}-z  } x_{2,j-1}^n$
on the rectangular
$[x_{i-1},x_i) \times I_{2,j-1}$, $1\leq i\leq N$, $1\leq j\leq M$.
Then
$$ D_a = \lim_{d\to 0} \lim_{r\to 0}
\sum_{i=1}^N \sum_{j=1}^M
\frac{ x_{1,i-1}^m }{ x_{1,i-1}-z } x_{2,j-1}^n
\left( \mathbf{E}_k ([x_{1,i-1},x_{1,i})) F(I_{2,j}) x_{0,0}, x_{0,0} \right)_H $$
$$ =
\lim_{d\to 0} \lim_{r\to 0}
\left( \sum_{i=1}^N
\frac{ x_{1,i-1}^m }{ x_{1,i-1}-z }  \mathbf{E}_k ([x_{i-1},x_i))
\sum_{j=1}^M
x_{2,j}^n F(I_{2,j}) x_{0,0}, x_{0,0} \right)_H $$
$$ =
\left( \int_{[-a,a)}
\frac{ x_1^m }{ x_1-z }  d\mathbf{E}_k \int_{[-R,R]}
x_2^n dF x_{0,0}, x_{0,0} \right)_H. $$
Let $n = n_1+n_2$, $n_1,n_2\in \mathbb{Z}_+$. Then we may write:
$$ D_a = \left( B^{n_1} \int_{[-a,a)}
\frac{ x_1^m }{ x_1-z }  d\mathbf{E}_k B^{n_2} x_{0,0}, x_{0,0} \right)_H $$
$$ =
\left( \int_{[-a,a)}
\frac{ x_1^m }{ x_1-z }  d\mathbf{E}_k x_{0,n_2}, x_{0,n_1} \right)_H. $$
By~(\ref{f3_28}) we get
$$ \int_\Pi \frac{x_1^m}{x_1-z} x_2^n d\mu =
\lim_{a\to +\infty} D_a =
\lim_{a\to +\infty}\left( \int_{[-a,a)}
\frac{ x_1^m }{ x_1-z }  d \widehat{E}_k x_{0,n_2}, x_{0,n_1} \right)_{\widehat H_k} $$
$$ = \left( \int_\mathbb{R}
\frac{ x_1^m }{ x_1-z }  d\widehat{E}_k x_{0,n_2}, x_{0,n_1} \right)_{\widehat H_k}
= \left( \widehat{A}_k^{m_2} R_z(\widehat{A}_k) \widehat{A}_k^{m_1} x_{0,n_2}, x_{0,n_1} \right)_{\widehat H_k}
$$
\begin{equation}
\label{f3_29}
= \left( R_z(\widehat{A}_k) x_{m_1,n_2}, x_{m_2,n_1} \right)_H,
\end{equation}
where $m_1,m_2\in \mathbb{Z}_+:\ m_1+m_2 = m$,
and $\widehat A_k$ is a self-adjoint extension of $A$ in a Hilbert space $\widehat H_k\supseteq H$ such that
its spectral measure $\widehat E_k$ generates $\mathbf{E}_k$: $\mathbf{E}_k = P^{\widehat H_k}_H \widehat E_k$;
$k=1,2$.

\noindent
Relation~(\ref{f3_29}) shows that the generalized resolvents corresponding to $\mathbf{E}_k$, $k=1,2$, coincide.
This means that the spectral measures $\mathbf{E}_1$ and $\mathbf{E}_2$ coincide. We obtained a contradiction.
This completes the proof.
$\Box$

\begin{dfn}
\label{d3_1}
A solution $\mu$ of the moment problem~(\ref{f1_1}) is said to be {\bf canonical}
if it is generated by relation~(\ref{f3_4}) where $\mathbf{E}$ is an {\bf orthogonal}
spectral measure of $A$ which commutes with the spectral measure of $B$. Orthogonal spectral measures
are those measures which are the spectral measures of self-adjoint extensions of $A$ inside $H$.
\end{dfn}
Let the moment problem~(\ref{f1_1}) be given and conditions~(\ref{f2_1}),(\ref{f2_2}) hold.
Let us describe canonical solutions of the two-dimensional moment problem in a strip.
In the proof of Theorem~\ref{t2_2} we have constructed one canonical solution, see relation~(\ref{f2_31}).
Let $\mu$ be an arbitrary canonical solution and $\mathbf{E}$ be the corresponding orthogonal spectral
measure of $A$. Let $\widetilde A$ be the self-adjoint operator in $H$ which corresponds to $\mathbf{E}$.
Consider the Cayley transformation of $\widetilde A$:
\begin{equation}
\label{f3_30}
U_{\widetilde A} = (\widetilde A + iE_H)(\widetilde A - iE_H)^{-1} \supseteq V_A,
\end{equation}
where $V_A$ is defined by~(\ref{f2_16}).
Since $\mathbf{E}$ commutes with the spectral measure $F$ of $B$, then $U_{\widetilde A}$ commutes
with $B$ and with $U_B$.
By relation~(\ref{f2_22}) the operator $U_{\widetilde A}$ have the following form:
\begin{equation}
\label{f3_31}
U_{\widetilde A} = V_A \oplus \widetilde U_{2,4},
\end{equation}
where $\widetilde U_{2,4}$ is an isometric operator which maps $H_2$ onto $H_4$, and
commutes with $U_B$.
Let the operator $U_{2,4}$ be defined by~(\ref{f2_26}). Then the following operator
\begin{equation}
\label{f3_32}
U_2 = U_{2,4}^{-1} \widetilde U_{2,4},
\end{equation}
is a unitary operator in $H_2$ which commutes with $U_{B;H_2}$.

Denote by $\mathbf{S}(U_B;H_2)$ a set of all unitary operators in $H_2$ which commute with $U_{B;H_2}$.
Choose an arbitrary operator $\widehat U_2\in \mathbf{S}(U_B;H_2)$. Define
$\widehat U_{2,4}$ by the following relation:
\begin{equation}
\label{f3_33}
\widehat U_{2,4} = U_{2,4} \widehat U_2.
\end{equation}
Notice that $\widehat U_{2,4} U_B h = \widehat U_{2,4} U_B h$, $h\in H_2$.
Then we define a unitary operator $U = V_A \oplus \widehat U_{2,4}$ and its Cayley transformation
$\widehat A$ which commute with the operator $B$.
Repeating arguments before~(\ref{f2_31}) we get a canonical solution of the moment problem.

\noindent
Thus, all canonical solutions of the Devinatz moment problem are generated by operators
$\widehat U_2\in \mathbf{S}(U_B;H_2)$. Notice that different operators $U',U''\in \mathbf{S}(U_B;H_2)$ produce different orthogonal spectral measures $\mathbf{E}',\mathbf{E}$. By Theorem~\ref{t3_1},
these spectral measures produce different solutions of the moment problem.

Recall some definitions from~\cite{cit_9000_BS}.
A pair $(Y,\mathfrak{A})$, where $Y$ is an arbitrary set and $\mathfrak{A}$ is a fixed
$\sigma$-algebra of subsets of $A$ is said to be a {\it measurable space}.
A triple $(Y,\mathfrak{A},\mu)$, where $(Y,\mathfrak{A})$ is a measurable space and
$\mu$ is a measure on $\mathfrak{A}$ is said to be a {\it space with a measure}.

Let $(Y,\mathfrak{A})$ be  a measurable space, $\mathbf{H}$ be a Hilbert space and
$\mathcal{P}=\mathcal{P}(\mathbf{H})$ be a set of all orthogonal projectors in $\mathbf{H}$.
A countably additive mapping $E:\ \mathfrak{A}\rightarrow \mathcal{P}$, $E(Y) = E_{\mathbf{H}}$,
is said to be a {\it spectral measure} in $\mathbf{H}$.
A set $(Y,\mathfrak{A},H,E)$ is said to be a {\it space with a spectral measure}.
By $S(Y,E)$ one means a set of all $E$-measurable $E$-a.e. finite complex-valued functions on $Y$.

Let $(Y,\mathfrak{A},\mu)$ be a separable space with a $\sigma$-finite measure and
to $\mu$-everyone $y\in Y$ it corresponds a Hilbert space $G(y)$. A function
$N(y) = \dim G(y)$ is called the {\it dimension function}.
It is supposed to be $\mu$-measurable. Let $\Omega$ be a set of vector-valued functions $g(y)$ with
values in $G(y)$ which are defined $\mu$-everywhere and are measurable with respect to some
base of measurability. A set of (classes of equivalence) of such functions with
the finite norm
\begin{equation}
\label{f3_34}
\| g \|^2_{\mathcal{H}} = \int |g(y)|^2_{G(y)} d\mu(y) <\infty
\end{equation}
form a Hilbert space $\mathcal{H}$ with the scalar product given by
\begin{equation}
\label{f3_35}
( g_1,g_2 )_{\mathcal{H}} = \int (g_1,g_2)_{G(y)} d\mu(y).
\end{equation}
The space $\mathcal{H}= \mathcal{H}_{\mu,N} = \int_Y \oplus G(y) d\mu(y)$
is said to be a {\it direct integral of Hilbert spaces}.
Consider the following operator
\begin{equation}
\label{f3_36}
\mathbf{X}(\delta) g = \chi_\delta g,\qquad g\in \mathcal{H},\ \delta\in \mathfrak{A},
\end{equation}
where $\chi_\delta$ is the characteristic function of the set $\delta$.
The operator $\mathbf{X}$ is a spectral measure in $\mathcal{H}$.

Let $t(y)$ be a measurable operator-valued function with values in $\mathbf{B}(G(y))$ which is
$\mu$-a.e. defined and $\mu-\sup \|t(y)\|_{G(y)} < \infty$. The operator
\begin{equation}
\label{f3_37}
T:\ g(y) \mapsto t(y)g(y),
\end{equation}
is said to be {\it decomposable}. It is a bounded operator in $\mathcal{H}$ which commutes with
$\mathbf{X}(\delta)$, $\forall\delta\in \mathfrak{A}$.
Moreover, every bounded operator in $\mathcal{H}$ which commutes with
$\mathbf{X}(\delta)$, $\forall\delta\in \mathfrak{A}$, is decomposable~\cite{cit_9000_BS}.
In the case $t(y) = \varphi(y)E_{G(y)}$, where $\varphi\in S(Y,\mu)$, we set $T =: Q_\varphi$.
The decomposable operator is unitary if and only if $\mu$-a.e. the operator $t(y)$ is unitary.

Return to the investigation of canonical solutions. Consider the spectral measure
$F_2$ of the operator $U_{B;H_2}$ in $H_2$.
There exists an element $h\in H_2$ of the maximal type, i.e. the non-negative Borel measure
\begin{equation}
\label{f3_38}
\mu(\delta) := (F_2(\delta)h,h),\qquad \delta\in \mathfrak{B}([-\pi,\pi)),
\end{equation}
has the maximal type between all such measures (generated by other elements of $H_2$). This type
is said to be the {\it spectral type} of the measure $F_2$.
Let $N_2$ be the multiplicity function of the measure $F_2$. Then there exists a unitary transformation $W$
of the space $H_2$ on $\mathcal{H}=\mathcal{H}_{\mu,N_2}$ such that
\begin{equation}
\label{f3_39}
W U_{B;H_2} W^{-1} = Q_{e^{iy}},\qquad W F_2(\delta) W^{-1} = \mathbf{X}(\delta).
\end{equation}
Notice that $\widehat U_2\in \mathbf{S}(U_B;H_2)$ if and only if
the operator
\begin{equation}
\label{f3_40}
V_2 := W \widehat U_2 W^{-1},
\end{equation}
is unitary and commutes with $\mathbf{X}(\delta)$, $\forall\delta\in \mathfrak{[-\pi,\pi)}$.
The latter is equivalent to the condition that $V_2$ is decomposable and the values of the corresponding
operator-valued function $t(y)$ are $\mu$-a.e. unitary operators.
A set of all decomposable operators in $\mathcal{H}$ such that the values of the corresponding
operator-valued function $t(y)$ are $\mu$-a.e. unitary operators we denote by $\mathbf{D}(U_B;H_2)$.

\begin{thm}
\label{t3_2}
Let the moment problem~(\ref{f1_1}) be given. In the conditions of Theorem~\ref{t3_1} all
canonical solutions of the moment problem have the  form~(\ref{f3_4}) where the spectral
measures $\mathbf{E}$ of the operator $A$ are constructed by operators from $\mathbf{D}(U_B;H_2)$.
Namely, for an arbitrary $V_2\in \mathbf{D}(U_B;H_2)$ we set $U_2 = W^{-1} V_2 W$,
$\widehat U_{2,4} = U_{2,4} \widehat U_2$, $U = V_A \oplus \widehat U_{2,4}$,
$\widehat A = i(U+E_H)(U-E_H)^{-1}$, and then $\mathbf{E}$ is the spectral measure of $\widehat A$.

\noindent
Moreover, the correspondence between $\mathbf{D}(U_B;H_2)$ and a set of all canonical solutions of
the moment problem is bijective.
\end{thm}
{\bf Proof. } The proof follows from the previous considerations.
$\Box$

Consider the moment problem~(\ref{f1_1}) and suppose that conditions~(\ref{f2_1}),(\ref{f2_2}) hold.
Let us turn to a parameterization of all solutions of the moment problem.
We shall use Theorem~\ref{t3_1}. Consider relation~(\ref{f3_4}). The spectral measure $\mathbf{E}$
commutes with the operator $U_B$.
Choose an arbitrary $z\in \mathbb{C}\backslash \mathbb{R}$.
By virtue of relation~(\ref{f3_3}) we may write:
$$ (U_B\mathbf{R}_z(A) x,y)_H = (\mathbf{R}_z(A) x,U_B^*y)_H =
\int_{\mathbb{R}} \frac{1}{t-z} d(\mathbf{E}(t) x,U_B^*y)_H $$
\begin{equation}
\label{f3_41}
= \int_{\mathbb{R}} \frac{1}{t-z} d(U_B\mathbf{E}(t) x,y)_H
= \int_{\mathbb{R}} \frac{1}{t-z} d(\mathbf{E}(t)U_B x,y)_H,\qquad x,y\in H;
\end{equation}
\begin{equation}
\label{f3_42}
(\mathbf{R}_z(A) U_Bx,y)_H = \int_{\mathbb{R}} \frac{1}{t-z} d(\mathbf{E}(t) U_Bx,y)_H,\qquad x,y\in H,
\end{equation}
where $\mathbf{R}_z(A)$ is  the generalized resolvent which corresponds to $\mathbf{E}$.
Therefore we get
\begin{equation}
\label{f3_43}
\mathbf{R}_z(A) U_B = U_B \mathbf{R}_z(A),\qquad z\in \mathbb{C}\backslash \mathbb{R}.
\end{equation}
On the other hand, if relation~(\ref{f3_43}) holds, then
\begin{equation}
\label{f3_44}
\int_{\mathbb{R}} \frac{1}{t-z} d(\mathbf{E} U_Bx,y)_H =
\int_{\mathbb{R}} \frac{1}{t-z} d(U_B\mathbf{E} x,y)_H,\quad x,y\in H,\ z\in \mathbb{C}\backslash \mathbb{R}.
\end{equation}
By the Stieltjes inversion formula~\cite{cit_10000_ST}, we obtain that $\mathbf{E}$ commutes with $U_B$.

\noindent
We denote by $\mathbf{M}(A,B)$ a set of all generalized resolvents $\mathbf{R}_z(A)$ of $A$ which satisfy
relation~(\ref{f3_43}).

Recall some known facts from~\cite{cit_4000_S} which we shall need here.
Let $K$ be a closed symmetric operator in a Hilbert space $\mathbf{H}$, with the domain $D(K)$,
$\overline{D(K)} = \mathbf{H}$.  Set $N_\lambda = N_\lambda(K) = \mathbf{H}
\ominus \Delta_K(\lambda)$, $\lambda\in \mathbb{C}\backslash \mathbb{R}$.

Consider an
arbitrary bounded linear operator $C$, which maps $N_i$ into $N_{-i}$.
For
\begin{equation}
\label{f3_45}
g = f + C\psi - \psi,\qquad f\in D(K),\ \psi\in N_i,
\end{equation}
we set
\begin{equation}
\label{f3_46}
K_C g = Kf + i C \psi + i \psi.
\end{equation}
The operator $K_C$ is said to be a {\it quasiself-adjoint extension of the operator $K$, defined by
the operator $K$}.

The following theorem can be found in~\cite[Theorem 7]{cit_4000_S}:
\begin{thm}
\label{t3_3}
Let $K$ be a closed symmetric operator in a Hilbert space $\mathbf{H}$ with the domain $D(K)$,
$\overline{D(K)} = \mathbf{H}$.
All generalized resolvents of the operator $K$ have the following form:
\begin{equation}
\label{f3_47}
\mathbf R_\lambda (K) = \left\{ \begin{array}{cc} (K_{F(\lambda)} - \lambda E_\mathbf{H})^{-1}, &
\mathop{\rm Im}\nolimits\lambda > 0\\
(K_{F^*(\overline{\lambda}) } - \lambda E_\mathbf{H})^{-1}, & \mathop{\rm Im}\nolimits\lambda < 0 \end{array}\right.,
\end{equation}
where $F(\lambda)$ is an analytic in $\mathbb{C}_+$ operator-valued function, which values are contractions
which map $N_i(A) = H_2$ into $N_{-i}(A) = H_4$ ($\| F(\lambda) \|\leq 1$),
and $K_{F(\lambda)}$ is the quasiself-adjoint extension of $K$ defined by $F(\lambda)$.

On the other hand, for any operator function $F(\lambda)$ having the above properties there corresponds by
relation~(\ref{f3_47}) a generalized resolvent of $K$.
\end{thm}
Observe that the correspondence between all generalized resolvents and functions $F(\lambda)$ in
Theorem~\ref{t3_3} is bijective~\cite{cit_4000_S}.

Return to the study of the moment problem~(\ref{f1_1}).
Let us describe the set $\mathbf{M}(A,B)$. Choose an arbitrary $\mathbf{R}_\lambda\in \mathbf{M}(A,B)$.
By~(\ref{f3_47}) we get
\begin{equation}
\label{f3_48}
\mathbf{R}_\lambda = (A_{F(\lambda)} - \lambda E_H)^{-1},\qquad \mathop{\rm Im}\nolimits\lambda > 0,
\end{equation}
where $F(\lambda)$ is an analytic in $\mathbb{C}_+$ operator-valued function, which values are contractions
which map $H_2$ into $H_4$, and
$A_{F(\lambda)}$ is the quasiself-adjoint extension of $A$ defined by $F(\lambda)$.
Then
$$ A_{F(\lambda)} = \mathbf{R}_\lambda^{-1} + \lambda E_H,\qquad \mathop{\rm Im}\nolimits\lambda > 0. $$
By virtue of relation~(\ref{f3_43}) we obtain
\begin{equation}
\label{f3_49}
U_BA_{F(\lambda)} h = A_{F(\lambda)} U_B h,\qquad h\in D(A_{F(\lambda)}),\ \lambda\in \mathbb{C}_+.
\end{equation}
Consider the following operators
\begin{equation}
\label{f3_50}
W_{\lambda} := (A_{F(\lambda)} + iE_H)(A_{F(\lambda)} - iE_H)^{-1} =
E_H + 2i(A_{F(\lambda)} - iE_H)^{-1},
\end{equation}
\begin{equation}
\label{f3_51}
V_A = (A +iE_H)(A - iE_H)^{-1} =
E_H + 2i(A - iE_H)^{-1},
\end{equation}
where $\lambda\in \mathbb{C}_+$.
Notice that (\cite{cit_4000_S})
\begin{equation}
\label{f3_52}
W_{\lambda} = V_A \oplus F(\lambda),\qquad \lambda\in \mathbb{C}_+.
\end{equation}
The operator $(A_{F(\lambda)} - iE_H)^{-1}$ is defined
on the whole $H$, see~\cite[p.79]{cit_4000_S}.
By relation~(\ref{f3_49}) we obtain
\begin{equation}
\label{f3_53}
U_B (A_{F(\lambda)} - iE_H)^{-1} h =
(A_{F(\lambda)} - iE_H)^{-1} U_B h,\qquad h\in H,\ \lambda\in \mathbb{C}_+.
\end{equation}
Then
\begin{equation}
\label{f3_54}
U_B W_\lambda = W_\lambda U_B,\qquad \lambda\in \mathbb{C}_+.
\end{equation}
Recall that by Proposition~\ref{p2_1} the operator $U_B$ reduces the subspaces $H_j$, $1\leq j\leq 4$,
and $U_BV_A = V_A U_B$. If we choose an arbitrary $h\in H_2$ and apply relations~(\ref{f3_54}),(\ref{f3_52}),
we get
\begin{equation}
\label{f3_55}
U_B F(\lambda) = F(\lambda) U_B,\qquad \lambda\in \mathbb{C}_+.
\end{equation}
Denote by $\mathbf{F}(A,B)$ a set of all analytic in $\mathbb{C}_+$ operator-valued functions
which values are contractions which map $H_2$ into $H_4$ and which satisfy relation~(\ref{f3_55}).
Thus, for an arbitrary $\mathbf{R}_\lambda\in \mathbf{M}(A,B)$ the corresponding function
$F(\lambda)$ belongs to $\mathbf{F}(A,B)$.

\noindent
On the other hand, choose an arbitrary $F(\lambda)\in \mathbf{F}(A,B)$.
Then we derive~(\ref{f3_54}) with $W_\lambda$ defined by~(\ref{f3_50}). Then we get~(\ref{f3_53}),(\ref{f3_49})
and therefore
\begin{equation}
\label{f3_56}
U_B \mathbf{R}_\lambda  = \mathbf{R}_\lambda U_B,\qquad \lambda\in \mathbb{C}_+.
\end{equation}
Calculating the conjugate operators for the both sides of the last equality we conclude that this
relation holds for all $\lambda\in \mathbb{C}$.

\noindent
Consider the spectral measure $F_2$ of the operator $U_{B;H_2}$ in $H_2$. We shall use
relation~(\ref{f3_39}).
Observe that $F(\lambda)\in \mathbf{F}(A,B)$ if and only if
the operator-valued function
\begin{equation}
\label{f3_57}
G(\lambda) := W U_{2,4}^{-1} F(\lambda) W^{-1},\qquad \lambda\in \mathbb{C}_+,
\end{equation}
is analytic in $\mathbb{C}_+$ and has values which are
contractions in $\mathcal{H}$ which commute with $\mathbf{X}(\delta)$, $\forall\delta\in \mathfrak{[-\pi,\pi)}$.

This means that for an arbitrary $\lambda\in \mathbb{C}_+$ the operator
$G(\lambda)$ is decomposable and the values of the corresponding
operator-valued function $t(y)$ are $\mu$-a.e. contractions.
A set of all decomposable operators in $\mathcal{H}$ such that the values of the corresponding
operator-valued function $t(y)$ are $\mu$-a.e. contractions we denote by $\mathrm{T}(B;H_2)$.
A set of all analytic in $\mathbb{C}_+$ operator-valued functions $G(\lambda)$ with values
in $\mathrm{T}(B;H_2)$ we denote by $\mathbf{G}(A,B)$.

\begin{thm}
\label{t3_4}
Let the two-dimensional moment problem in a strip~(\ref{f1_1}) be given. In the conditions of Theorem~\ref{t3_1} all
solutions of the moment problem have the  form~(\ref{f3_4}) where the spectral
measures $\mathbf{E}$ of the operator $A$ are defined by the corresponding generalized
resolvents $\mathbf{R}_\lambda$ which are constructed by the following relation:
\begin{equation}
\label{f3_58}
\mathbf{R}_\lambda = (A_{F(\lambda)} - \lambda E_H)^{-1},\qquad \mathop{\rm Im}\nolimits\lambda > 0,
\end{equation}
where $F(\lambda) = U_{2,4} W^{-1} G(\lambda) W$, $G(\lambda)\in \mathbf{G}(A,B)$.

\noindent
Moreover, the correspondence between $\mathbf{G}(A,B)$ and a set of all solutions of
the moment problem is bijective.
\end{thm}
{\bf Proof. } The proof follows from the previous considerations.
$\Box$

\section{The complex moment problem in a strip.}
In this Section we shall analyze the following problem:
to find a non-negative Borel measure $\sigma$ in a strip
$$ \Psi = \Psi(R) = \{ z\in \mathbb{C}:\ | \mathop{\rm Im}\nolimits z | \leq R \},\qquad R>0, $$
such that
\begin{equation}
\label{f4_1}
\int_\Psi z^m \overline{z}^n d\sigma = a_{m,n},\qquad m,n\in \mathbb{Z}_+,
\end{equation}
where $\{ a_{m,n} \}_{m,n\in \mathbb{Z}_+}$ is a prescribed sequence of complex numbers.
This problem is said to be {\bf the complex moment problem in a strip}.
Of course, the strips $\Psi$ and $\Pi$ are the same sets in accordance with the canonical identification
of $\mathbb{C}$ with $\mathbb{R}^2$:
\begin{equation}
\label{f4_2}
z = x_1 + x_2 i,\ x_1 = \mathop{\rm Re}\nolimits z,\ x_2 = \mathop{\rm Im}\nolimits z,\quad
z\in \mathbb{C},\ (x_1,x_2)\in \mathbb{R}^2.
\end{equation}
Let $\sigma$ be a solution of the complex moment problem in a strip~(\ref{f4_1}).
The measure $\sigma$, viewed as a measure in $\mathbb{R}^2$, we shall denote by $\mu_\sigma$. Then
$$ s_{m,n} := \int_\Pi x_1^m x_2^n d\mu_\sigma = \int_\Psi \left( \frac{z+\overline{z}}{2} \right)^m
\left( \frac{z-\overline{z}}{2i} \right)^n d\sigma $$
$$ = \frac{1}{2^m (2i)^{n}} \sum_{k=0}^m \sum_{j=0}^n
C^m_k C^n_j (-1)^{n-j}\int_\Psi z^{k+j} \overline{z}^{m-k+n-j} d\sigma $$
\begin{equation}
\label{f4_3}
= \frac{1}{2^m (2i)^{n}} \sum_{k=0}^m \sum_{j=0}^n (-1)^{n-j} C^m_k C^n_j a_{k+j,m-k+n-j},
\end{equation}
where $C^n_k = \frac{n!}{k!(n-k)!}$.
Then
$$ a_{m,n} = \int_\Psi z^m \overline{z}^n d\sigma = \int_\Pi (x_1+ix_2)^m (x_1-ix_2)^n d\mu_\sigma  $$
$$ = \sum_{r=0}^m \sum_{l=0}^n C^m_r C^n_l (-1)^{n-l} \int_\Pi x_1^{r+l} (ix_2)^{m-r+n-l} d\mu_\sigma $$
$$ = \sum_{r=0}^m \sum_{l=0}^n C^m_r C^n_l (-1)^{n-l} i^{m-r+n-l}
s_{r+l,m-r+n-l}; $$
and therefore
\begin{equation}
\label{f4_4}
a_{m,n} = \sum_{r=0}^m \sum_{l=0}^n C^m_r C^n_l (-1)^{n-l} i^{m-r+n-l}
s_{r+l,m-r+n-l},\quad m,n\in \mathbb{Z}_+,
\end{equation}
where
\begin{equation}
\label{f4_5}
s_{m,n} = \frac{1}{2^{m+n}} \sum_{k=0}^m \sum_{j=0}^n (-1)^{j} C^m_k C^n_j a_{k+j,m-k+n-j},\quad m,n\in \mathbb{Z}_+.
\end{equation}
Since $\mu_\sigma$ is a solution of the two-dimensional moment problem in a strip, then
conditions~(\ref{f2_1}),(\ref{f2_2}) hold.
\begin{thm}
\label{t4_1}
Let the complex moment problem in a strip~(\ref{f4_1}) be given.
This problem has a solution if an only if conditions~(\ref{f2_1}),(\ref{f2_2}) and~(\ref{f4_4})
with $s_{m,n}$ defined by~(\ref{f4_5})  hold
for arbitrary complex numbers $\alpha_{m,n}$ such that all but finite numbers are zeros.
\end{thm}
{\bf Proof. }
It remains to prove the sufficiency. Suppose that for the complex moment problem in a strip~(\ref{f1_1})
conditions~(\ref{f2_1}),(\ref{f2_2}) and~(\ref{f4_4})  hold. By Theorem~\ref{t2_2} we obtain that
there exists a solution $\mu$ of the two-dimensional moment problem with moments $s_{m,n}$ defined
by~(\ref{f4_5}).
The measure $\mu$, viewed as a measure in $\mathbb{C}$, we shall denote by $\sigma_\mu$. Then
$$ \int_\Psi z^m \overline{z}^n d\sigma_\mu = \int_\Pi (x_1+ix_2)^m (x_1-ix_2)^n d\mu  $$
$$ = \sum_{r=0}^m \sum_{l=0}^n C^m_r C^n_l (-1)^{n-l} \int_\Pi x_1^{r+l} (ix_2)^{m-r+n-l} d\mu $$
$$ = \sum_{r=0}^m \sum_{l=0}^n C^m_r C^n_l (-1)^{n-l} i^{m-r+n-l}
s_{r+l,m-r+n-l} = a_{m,n}, $$
where the last equality follows from~(\ref{f4_4}).
$\Box$

\begin{thm}
\label{t4_2}
Let the complex moment problem in a strip~(\ref{f4_1}) be given and
conditions~(\ref{f2_1}),(\ref{f2_2}) and~(\ref{f4_4})  hold for arbitrary
complex numbers $\alpha_{m,n}$ such that all but finite numbers are zeros.
All solutions of the moment problem~(\ref{f4_1}) are solutions of the two-dimensional moment problem~(\ref{f1_1})
with $s_{m,n}$ defined by~(\ref{f4_5}), viewed as measures on $\mathbb{C}$.
Therefore all solutions of the moment problem~(\ref{f4_1}) have a parameterization
provided by Theorem~\ref{t3_4}.
\end{thm}
{\bf Proof. } The proof is straightforward.
$\Box$

\vspace{1.5cm}

Sergey M. Zagorodnyuk

School of Mathematics and Mekhanics

Karazin Kharkiv National University

Kharkiv, 61077

Ukraine

Sergey.M.Zagorodnyuk@univer.kharkov.ua

\begin{center}
\bf
The two-dimensional moment problem in a strip.
\end{center}

\begin{center}
\bf
S.M. Zagorodnyuk
\end{center}

In this paper we study the two-dimensional moment problem in a strip
$\Pi(R) = \{ (x_1,x_2)\in \mathbb{R}^2:\ |x_2| \leq R \}$, $R>0$.
We obtained a solvability criterion for this moment problem.
We derived a parameterization of all solutions of the moment problem.
An abstract operator approach and results of Godi\v{c}, Lucenko and
Shtraus are used.

\vspace{1cm}
Key words: moment problem, measure, generalized resolvent.
\vspace{1cm}

MSC 2000: 44A60, 30E05.

}
\end{document}